\documentclass[11pt]{article}

\usepackage{amsmath}
\usepackage{amssymb}
\usepackage{amsthm}
\usepackage{amsrefs}
\usepackage{mathtools}
\usepackage{mathrsfs}
\usepackage{tikz}
\usepackage[margin=1in]{geometry}
\usepackage{setspace}
\usepackage{xcolor}
\usepackage{ytableau}

\newtheorem{theorem}{Theorem}

\newtheorem{proposition}{Proposition}

\theoremstyle{definition}

\newtheorem{conjecture}{Conjecture}
\newtheorem{remark}{Remark}

\title{An Overpartition Companion of Andrews and Keith's 2-colored $q$-series Identity}
\author{Hunter Waldron}
\date{}

\newcommand\qbin[3]{{\genfrac{[}{]}{0pt}{1}{#1}{#2}}_{\displaystyle #3}}

\begin{document}

\maketitle
  
\begin{abstract}
Andrews and Keith recently produced a general Schmidt type partition theorem using a novel interpretation of Stockhofe's bijection, which
they used to find new $q$-series identities. This includes an identity for a trivariate 2-colored partition generating function.
In this paper, their Schmidt type theorem is further generalized akin to how Franklin classically extended
Glaisher's theorem. As a consequence, we obtain a companion to Andrews and Keith's 2-colored identity for overpartitions.
These identities appear to be special cases of a much more general result. 
\end{abstract}

\section{Schmidt Type and Colored Partitions}

\indent There has been a growing interest
\cite{AndrewsKeith, WangXu, Alladi2, WangZhangZhou, ChernYee, Ji, BerkovichUncu, LiYee, BridgesUncu} following the expository work of
Andrews and Paule \cite{AndrewsPaule}
in researching the combinatorics of Schmidt type partitions, which are partitions 
weighted by the sum of parts having only indices from some strict subset of the positive integers.
Their namesake Schmidt was the first to observe Theorem \ref{Schmidts_theorem}, which has since witnessed
numerous proofs (notably, Mork produced a bijection \cite{Mork} using diagonal hooks, similar to Bessenrodt's variant
of Sylvester's classic bijection \cite{Bessenrodt}). 

\begin{theorem}[Schmidt \cite{Schmidt}]
For all $n \geq 0$, the number of partitions $\lambda$ with distinct parts
such that $\lambda_1 + 
\lambda_3 + \lambda_5 + \cdots = n$ is equal to the number of partitions of $n$. Equivalently,
\begin{equation*}
\smashoperator{\sum_{\lambda \in \mathcal{D}}} q^{\displaystyle \lambda_1 + \lambda_3 + \lambda_5 + \cdots}
= \smashoperator{\sum_{\lambda \in \mathcal{P}}} q^{\displaystyle \lvert \lambda \rvert}.
\end{equation*}
\label{Schmidts_theorem}
\end{theorem}

The sets $\mathcal{P}$ and $\mathcal{D}$ shown above are the set of all partitions and the set of partitions with distinct parts, respectively.
For Schmidt weights to be defined,
we must take any partition $\lambda \in \mathcal{P}$ bearing them to be a weakly decreasing sequence
$\lambda = (\lambda_1, \lambda_2, \dots, \lambda_{\ell(\lambda)})$ of $\ell(\lambda)$
positive integers with $\lambda_k = 0$ whenever $k > \ell(\lambda)$, a convention used implicitly
throughout this paper. The parts of $\lambda$ are the positive $\lambda_k$, the size is the sum of
all parts $\lvert \lambda \rvert = \lambda_1 + \lambda_2 + \lambda_3 + \cdots$, and the length is
the number of parts $\ell(\lambda)$. When $\lambda$ has size $n$ we say $\lambda$ is a partition of $n$, and
there is a unique partition with no parts $\lambda = ()$, which is the only partition satisfying $\lvert \lambda \rvert = 0$ or $\ell(\lambda) = 0$.

\begin{remark}
Bowman and Alladi each found results before Schmidt that can be considered
equivalent to Theorem \ref{Schmidts_theorem} (see \cite{Bowman, Alladi1, Alladi2}).
\end{remark}

A result analogous to Theorem \ref{Schmidts_theorem} was found by Uncu \cite{Uncu} as a special
case of an identity using Boulet-Stanley weights \cite{Boulet}.
\begin{theorem}[Uncu \cite{Uncu}]
For all $n \geq 0$, the number of partitions $\lambda$ such that $\lambda_1 + \lambda_3 + \lambda_5 + \cdots = n$
is equal to the number of 2-colored partitions of $n$. Equivalently,
\begin{equation*}
\smashoperator{\sum_{\lambda \in \mathcal{P}}} q^{\displaystyle \lambda_1 + \lambda_3 + \lambda_5 + \cdots}
= \smashoperator{\sum_{\lambda \in \mathcal{P}_2}} q^{\displaystyle \lvert \lambda \rvert}.
\end{equation*}
\label{uncus_theorem}
\end{theorem}

Here we introduce the language of colored partitions, where a part size is said to appear
in $t$ colors, which we take to be some fixed set of positive integers $\{c_1, \dots, c_t\}$, if
distinct markings of all such parts using these colors as subscripts in
weakly decreasing order are distinguished. A $t$-colored partition may have any part
size appear in $t$ colors, and the set of all such partitions without restrictions is written $\mathcal{P}_t$.
For example, with colors $\{1, 2\}$ the 2-colored partitions of $3$ are
\begin{equation*}
\begin{gathered}
(3_2),\; (3_1),\; (2_2, 1_2),\; (2_2, 1_1),\; (2_1, 1_2),\; (2_2, 1_2),\;
(1_2, 1_2, 1_2), \\ (1_2, 1_2, 1_1), \;(1_2, 1_1, 1_1), \;\mathrm{and}\; (1_1, 1_1, 1_1).
\end{gathered}
\end{equation*}

By identifying $\mathcal{P}_1$ with $\mathcal{P}$, both Uncu and Schmidt's theorems hint 
that Schmidt weights may be related to colored partitions in general. Some more recent results show that they are indeed intimately linked.
In two papers by Li and Yee \cite{LiYee} and by Andrews and Keith \cite{AndrewsKeith},  strong generalizations of
Uncu's and Schmidt's theorems, respectively, were proven. Their results can be expressed in the following language.

\begin{theorem}[Li-Yee \cite{LiYee}]
For all $n, s, j \geq 1$ and $t \geq 2$, the number of partitions $\lambda \in \mathcal{P}$ such that
\begin{align*}
n &= \sum_{k \geq 0} \lambda_{tk+1} = \lambda_1 + \lambda_{t+1} + \lambda_{2t+1} + \cdots \\
\ell(\lambda) &= (s-1)t + j 
\end{align*}
is equal to  the number of $t$-colored partitions of $n$ where $s$ is the most number of times any color appears, and
$j$ is the largest such color.
\label{uncus_theorem_generalization}
\end{theorem}

\begin{theorem}[Andrews-Keith \cite{AndrewsKeith}]
Fix $m \geq 2$ and $S = \{s_1, \dots, s_i\} \subseteq \{1, 2, \dots, m\}$ with
$1 \in S$, and $S$ ordered $s_1 < s_2 < \cdots < s_i$.
For all $n \geq 1$, The partitions $\lambda$ with all part sizes appearing fewer than $m$ times such that
\begin{equation*}
\begin{aligned}
n &= \smashoperator{\sum_{\substack{k \equiv S \\ (\mathrm{mod}\; m)}}} \lambda_k =
\lambda_{s_1} + \cdots + \lambda_{s_i} + \lambda_{s_1+m} +
\cdots + \lambda_{s_i+m} + \cdots \\
\rho_j &= \smashoperator{\sum_{k \geq 0}} (\lambda_{mk+j} - \lambda_{mk+j+1}) =
\lambda_j - \lambda_{j+1} + \lambda_{j+m} - \lambda_{j+m+1} + \cdots
\end{aligned}
\end{equation*}
for $1 \leq j < m$ are equinumerous with the partitions of size $n$ in $\mathcal{P}$
where any parts congruent to  $k$ modulo $i$ appears in the $s_{k+1} - s_k$ colors
$\{s_k, \dots, s_{k+1} - 1\}$
where we take $s_{i+1} = m$, and parts of color $j$ appear $\rho_j$ times.
\label{AndrewsKeith_main}
\end{theorem}

\begin{remark}
The presentation of Theorem \ref{AndrewsKeith_main} shown here is slightly more general than the result that appears in \cite{AndrewsKeith}. Their proof
works just as well with edge cases that they omitted, when $m=2$ or $m \in S$.
\end{remark}

Results like Theorems \ref{uncus_theorem_generalization} and
\ref{AndrewsKeith_main} are examples of what are often called partition identities in the literature of partitions.
Although partition identities are interesting in their own right, they are often applied to find combinatorial proofs of $q$-series identities,
or even to derive new ones.
In \cite{AndrewsKeith}, Andrews and Keith did the latter, using some special cases of Theorem \ref{AndrewsKeith_main} to find
new $q$-series sum-product identities with the sum and product side corresponding to the Schmidt type and colored partitions, respectively.
Among their results is a trivariate identitiy that has a product side well-known to
be the generating function for 2-colored partitions
\begin{equation*}
\smashoperator{\sum_{\lambda \in \mathcal{P}_2}} t_1^{\ell_1(\lambda)} t_2^{\ell_2(\lambda)} q^{\displaystyle \lvert \lambda \rvert}
= \frac{1}{(t_1q; q)_\infty (t_2q; q)_\infty}
\end{equation*}
where $\ell_1(\lambda)$ and $\ell_2(\lambda)$ count the number of parts appearing in the respective colors.

\begin{theorem}[Andrews-Keith \cite{AndrewsKeith}] We have the equality
\begin{equation*}
\sum_{n \geq 0} \sum_{\substack{j + k \geq n \\ j, k \leq n}}
\frac{(-1)^{j+k+n}t_1^j t_2^k q^{{n \choose 2} + {j+1 \choose 2} + {k+1 \choose 2}}
\qbin{n}{n-j, n-k, j+k-n}{q}
}{(t_1q;q)_n (t_2q;q)_n (q;q)_n} =
\frac{1}{(t_1q; q)_\infty (t_2q; q)_\infty}.
\end{equation*}
\label{AndrewsKeith_example}
\end{theorem}

This is written as an equality of formal power series
using the standard notation for the $q$-Pochhammer symbol
$(z;q)_0 = 1$ and $(z;q)_n = \prod_{k=0}^{n-1}(1-zq^k)$ for $n>0$, with the limiting case as $n \to \infty$
written $(z;q)_\infty$. The
$q$-multinomial coefficient is given below on the left for $n = k_1 + \cdots + k_t$ with each $k_i \geq 0$ 
\begin{equation*}
\genfrac{[}{]}{0pt}{}{n}{k_1, \dots, k_t}_{\displaystyle q} = \frac{(q;q)_n}
{(q;q)_{k_1} \cdots (q;q)_{k_t}}\quad \quad
\genfrac{[}{]}{0pt}{}{n}{k}_{\displaystyle q} = \frac{(q;q)_n}{(q;q)_k(q;q)_{n-k}}
\end{equation*}
and the $q$-binomial coefficient, which appears later in this article, is the special case where $t=2$, and can unambiguously be written as shown
above on the right. 

\section{New Results}

The 2-colored partitions are related to the well-studied overpartitions $\overline{\mathcal{P}}$, introduced by Corteel and Lovejoy in \cite{CorteelLovejoy},
which are partitions in $\mathcal{P}$ where the first occurrence of any part size may be distinguished with an overline.
Overpartitions can be identified with 2-colored partitions where a fixed color may appear at most
once per part size. For example, the overpartitions of 3 are 

\begin{equation*}
(3),\; (\overline{3}), \;(2, 1), \;(\overline{2},1), \;(2, \overline{1}), \;(\overline{2},\overline{1}), \;(1,1,1), \;\mathrm{and}\;  (\overline{1},1,1).
\end{equation*}

Many well-known $q$-series identities have natural overpartition analogs (for a good example, Dousse and Kim found many overpartition analogs of $q$-binomial
identities \cite{DousseKim1, DousseKim2}). The main purpose of this paper is to provide an analog of Theorem
\ref{AndrewsKeith_example}, with the generating function for overpartitions
\begin{equation*}
\smashoperator{\sum_{\lambda \in \overline{\mathcal{P}}}} t_1^{\displaystyle \ell_o(\lambda)}
t_2^{\displaystyle \ell_n(\lambda)} q^{\displaystyle \lvert \lambda \rvert} = \frac{(-t_1q;q)_\infty}{(t_2q;q)_\infty}
\end{equation*}
appearing on the product side where $\ell_o(\lambda)$ and $\ell_{n}(\lambda)$ count respectively the number of overlined and non-overlined parts
in $\lambda$. Our main result is as follows.

\begin{theorem} We have the equality
\begin{equation*}
\sum_{n \geq 0} \sum_{\substack{j+k \geq n \\ j, k \leq n}}
\frac{ (-1)^{j+k+n}t_1^j t_2^k q^{{n \choose 2} + {k+1 \choose 2} + j^2-nj+j}
\qbin{n}{n-j, n-k, j+k-n}{q}
}{(t_2q;q)_n (q;q)_n} =
\frac{(-t_1q; q)_\infty}{(t_2q; q)_\infty}.
\end{equation*}
\label{overpartition_identity}
\end{theorem}

We ultimately prove Theorem \ref{overpartition_identity} using $q$-series techniques below in Section 3, but the identity was first conjectured only by establishing
an appropriate partition identity, which is most usefully written as the following equality of generating functions.
\begin{proposition}
Let  $\mathcal{D}_4$ be  the set of partitions with any part size repeating fewer than 4 times, and for any $\lambda 
\in \mathcal{D}_4$, let $e(\lambda)$ be the number of part sizes that appear 2 or 3 times. Then
\begin{equation*}
\sum_{\lambda \in \mathcal{D}_4} t_1^{\displaystyle e(\lambda)}
t_2^{\displaystyle \lambda_1 - \lambda_2 + \lambda_3
-\lambda_4 + \cdots} q^{\displaystyle \lambda_1 + \lambda_3 + \lambda_5 + \cdots} =
\smashoperator{\sum_{\lambda \in \overline{\mathcal{P}}}} t_1^{\displaystyle \ell_o(\lambda)}
t_2^{\displaystyle \ell_n(\lambda)} q^{\displaystyle \lvert \lambda \rvert}.
\end{equation*}
\label{overpartition_gf_identity}
\end{proposition}

\begin{remark}
This identity is indeed a generalization of Theorem \ref{Schmidts_theorem}, since only the partitions with distinct parts and with no overlines
contribute to the coefficient of $t_1^0$ on the left and right side, respectively.
\end{remark}

An example of the partitions being enumerated is provided here. The coefficient of $t_1t_2^2q^6$ is 6 since the relevant
Schmidt type partitions are
\begin{equation*}
( 5, 3, 1, 1),\; ( 4, 3, 1, 1, 1), \;( 3, 3, 3, 1), \;(4,4,2), \;( 3, 2, 2, 2, 1), \;\mathrm{and}\;  ( 4, 2, 2, 2),
\end{equation*}
and the  overpartitions are
\begin{equation*}
(\overline{4},1,1), \;(4, \overline{1}, 1), \;(\overline{2},2,2),
\;(\overline{3},2,1), \;(3,\overline{2},1), \;\mathrm{and} \;(3,2,\overline{1}).
\end{equation*}

Proposition \ref{overpartition_gf_identity} does not appear to be a special case of Theorem \ref{AndrewsKeith_main}, at least in the full trivariate generality.
A secondary purpose of this paper is to provide the following extension of Theorem \ref{AndrewsKeith_main} to allow arbitrary part repetition on the Schmidt
type side to be controlled instead of forbidden, from which partition identities like Proposition \ref{overpartition_gf_identity}
may be easily established as a corollary. In particular, Proposition \ref{overpartition_gf_identity} is the limiting case where we
take $S=\{1\}$, $m=2$, the $\alpha_i$ to be all positive integers, and set all $\beta_i = 2$.

\begin{theorem}
Fix $m$ and $S$ as in Theorem \ref{AndrewsKeith_main}. Let  $\alpha_1, \dots, \alpha_t$ be distinct positive integers and
$\beta_1, \dots \beta_t \geq 1$. Then the partitions $\lambda \in \mathcal{P}$ with the same conditions on the Schmidt weight $n$ and the $\rho_j$,
such that part sizes $\alpha_i$ appear fewer than $\beta_i m$ times and all others appear fewer than $m$ times, are equinumerous with the partitions
of size $n$ with the same conditions on the colors, such that in addition parts $i \alpha_i$ may appear in the color $m$ fewer than $\beta_i$ times.
\label{AndrewsKeith_extension}
\end{theorem}

\begin{remark}
The generalization of Theorem \ref{AndrewsKeith_main} to Theorem \ref{AndrewsKeith_extension} is similar to how Franklin
\cite{Franklin} classically extended Glaisher's theorem \cite{Glaisher}.
\end{remark}

We finish this section with a discussion of the future work that may draw from \cite{AndrewsKeith} and this paper.
All the $q$-series identities that Andrews and Keith proved in \cite{AndrewsKeith} came from quite special cases of their partition identity.
This is a strong reason to suspect that their results, along with Theorem \ref{overpartition_identity}, are the simplest cases of infinite
families of $q$-series identities, indexed by many parameters. As an example, we offer the following conjecture of a one parameter
infinite family, with Theorem \ref{AndrewsKeith_main} retrieved when $k=2$.

\begin{conjecture}
For each $k \geq 2$, there exists polynomials $f(n, i_1, \cdots, i_k, q)$ in $q$ with non-negative coefficients and a constant term 1 such that
\begin{align*}
\sum_{n \geq 0}\hspace{0.7cm} \smashoperator{\sum_{\substack{i_1 + \cdots + i_k \geq  n \\ 0 \leq i_1, \dots, i_k \leq n}}} \hspace{0.2cm}
&\frac{(-1)^{n+i_1+\cdots+i_k} t_1^{i_1}\cdots t_k^{i_k} q^{{n \choose 2} + {i_1 + 1 \choose 2} + \cdots + {i_k + 1 \choose 2}}f(n, i_1, \cdots, i_k, q)}
{(q;q)_n(t_1q;q)_n \cdots (t_kq;q)_n} \\ &= \frac{1}{(t_1q;q)_\infty( t_2q;q)_\infty\cdots ( t_k q; q)_\infty}.
\end{align*}
\end{conjecture}

This conjecture is supported by computational evidence, although the form these unknown polynomials might take in general is not clear.
Plainly they should be symmetric in the indices $i_1, \dots, i_k$, and for $k > 2$, if all but two of the $t_1, \dots, t_k$ are sent to 0, this identity
reduces to Theorem \ref{AndrewsKeith_main}. From this, we can at least say that $f(n, i_1, \dots, i_k, q)$ reduces to the $q$-multinomial coefficient when at most
two of these indices are non-zero. Some other examples are given here with $k=3$ and $n=4$:
\begin{align*}
 f(4, 2, 2, 2, q) &= (1 + q^2)  (1 + q + q^2)  (1 + 2 q + 5 q^2 + 5 q^3 + 5 q^4 + q^5), \\
 f(4, 1, 2, 3, q) &= (1 + q)  (1 + q^2)  (1 + q + q^2)  (1 + q + 2 q^2 + q^3), \\
 f(4, 2, 3, 4, q) &= (1 + q)  (1 + q^2)^2  (1 + q + q^2), \\
 f(4, 2, 2, 3, q) &= (1 + q)  (1 + q^2)  (1 + q + q^2)^2  (1 + 2  q^2).
\end{align*}

Naturally, Theorem \ref{overpartition_identity} may have a similar generalization, although what that might look like, and whether these polynomials
are involved is even less obvious.

The rest of the paper is organized into two additional sections. Directly below in Section 3, the proof of Theorem \ref{overpartition_identity} is given.
Then, in Section 4, a second proof of Theorem \ref{AndrewsKeith_main} is provided using generating functions,
and we provide a proof of Theorem \ref{AndrewsKeith_extension} by extending the bijection originally used by Andrews and Keith.

\section{Proof of Theorem \ref{overpartition_identity}}

\noindent \textit{Proof of Theorem \ref{overpartition_identity}}.
Using well-known generating function identities, we can write
\begin{equation*}
\frac{(-t_1q; q)_\infty}{(t_2q;q)_\infty} = \sum_{n \geq 0} \frac{t_1^nq^{n+1\choose 2}}{(q;q)_n}
\times \sum_{m \geq 0} \frac{t_2^m q^{m^2}}{(q;q)_m (t_2q;q)_m}
\end{equation*}
and so for any fixed $J \geq 0$, the coefficient of $t_1^J$ on the product side of Theorem \ref{overpartition_identity}  is 
\begin{equation*}
\frac{q^{J+1 \choose 2}}{(q;q)_J}\sum_{m \geq 0} \frac{t_2^m q^{m^2}}{(q;q)_m (qt_2;q)_m}.
\end{equation*}
We proceed now to show that the coefficient of $t_1^J$ is the same on the sum side.
Extracting the coefficient of $t_1^J$, we obtain
\begin{align*}
&(-1)^Jq^{J^2+J}\sum_{n \geq J} \sum_{\substack{J+k \geq n \\ k \leq n}}
\frac{ (-1)^{k+n}t_2^k q^{{n \choose 2} + {k+1 \choose 2} -nJ}
\qbin{n}{n-J, n-k, J+k-n}{q}
}{(t_2q;q)_n (q;q)_n} \\
=&q^{J+1 \choose 2}\sum_{n \geq 0} \sum_{k=n}^{n+J}
\frac{ (-1)^{k+n}t_2^k q^{{n \choose 2} + {k+1 \choose 2}}
\qbin{n+J}{n, n+J-k, k-n}{q}
}{(t_2q;q)_{n+J} (q;q)_{n+J}} \\
=&q^{J+1 \choose 2}\sum_{n \geq 0} \sum_{k=0}^{J}
\frac{ (-1)^{k}t_2^{n+k} q^{n^2 + {k +1 \choose 2} +nk}
\qbin{n+J}{n, J-k, k}{q}}{(t_2q;q)_{n+J} (q;q)_{n+J}} \\
=&\frac{q^{J+1 \choose 2}}{(q;q)_J} \sum_{n \geq 0} \frac{t_2^nq^{n^2}}{(t_2q;q)_{n+J} (q;q)_{n}}
\sum_{k=0}^J (-1)^kt_2^k q^{{k\choose 2} + (n+1)k }\qbin{J}{k}{q} 
\end{align*}
by shifting $n$ to $n+J$, $k$ to $n+k$, and then finally
breaking apart the $q$-multinomial coefficient and rewriting. From here we can invoke a standard result, 
Cauchy's $q$-binomial theorem (see \cite{Andrews}, Theorem 3.3)
\begin{equation*}
(z;q)_N = \sum_{k=0}^N (-1)^k z^k q^{k \choose 2} \qbin{N}{k}{q}
\end{equation*}
with $z=t_2q^{n+1}$ to write this in the desired form
\begin{align*}
&\frac{q^{J+1 \choose 2}}{(q;q)_J} \sum_{n \geq 0} \frac{t_2^nq^{n^2}}{(t_2q;q)_{n+J} (q;q)_{n}}
(t_2q^{n+1};q)_J \\
=&\frac{q^{J+1 \choose 2}}{(q;q)_J} \sum_{n \geq 0} \frac{t_2^nq^{n^2}}{(t_2q;q)_{n} (q;q)_{n}}.
\end{align*}
\qed

\section{On Theorem \ref{AndrewsKeith_main} and a Proof of the Extension.}

We start this section by introducing some notation as well as several maps between partitions that will be needed in the bijective proof of Theorem
\ref{AndrewsKeith_extension} below.

Any partition $\lambda \in \mathcal{P}$ may be represented with a Young diagram, which is a grid of
squares having $\lambda_i$ squares in the $i$th row down, with all rows vertically aligned on the left. The conjugate partition $\lambda'$ of
$\lambda$ is formed by taking $\lambda'_i$ to be the number of squares in the $i$th column from the left in
$\lambda$'s Young diagram (see Figure \ref{Young_diagram_example} for an example). Conjugation is a bijection from $\mathcal{P}$ to $\mathcal{P}$.

\begin{figure}[!ht]
\centering
\begin{tikzpicture}[scale=0.4]
\draw (0,0) grid (5, -1);
\draw (0,-1) grid (4, -2);
\draw (0,-2) grid (1, -6);
\draw (0,-2) grid (2,-3);
\draw (10,0) grid (16, -1);
\draw (10,-1) grid (13,-2);
\draw (10,-2) grid (12,-4);
\draw (10,-4) grid (11,-5);
\end{tikzpicture}
\caption{The Young diagram of $\lambda = (5,4,2,1,1,1)$ is shown on the left, with the conjugate $\lambda'=(6,3,2,2,1)$ on the right}
\label{Young_diagram_example}
\end{figure}

For $m \geq 2$, let $\mathcal{D}_m$ be the set of partitions with part sizes repeating
fewer than $m$ times, $\mathcal{R}_m$ be the set of partitions with no parts divisible by $m$, and $\mathcal{F}_m$ be the set of partitions
$\lambda$ such that $\lambda_i - \lambda_{i+1} < m$ for each $i \geq 1$. Equivalently, $\mathcal{F}_m = 
\{\lambda' \mid \lambda \in \mathcal{D}_m\}$.

Stockhofe gave a family of bijections indexed by $m \geq 2$ in \cite{Stockhofe}, which specialize to a map $\phi \colon \mathcal{F}_m \to \mathcal{R}_m$.
A definition of $\phi$ in English is provided in \cite{AndrewsKeith}, along with a detailed example, and so one is not provided here. What we need are
the facts that for any $\lambda \in \mathcal{D}_m$, $\lambda$ and  $\phi(\lambda')'$ will have the same Schmidt weight when the weight repeats
with period $m$, given by the set $S$ as in Theorem \ref{AndrewsKeith_main},
and $\phi(\lambda')'$ also preserves the statistic $\rho_j = \lambda_j - \lambda_{j+1} + \lambda_{j+m} - \lambda_{j+m+1} + \cdots$.

\begin{figure}[!ht]
\centering
\begin{tikzpicture}[scale=0.4]

\filldraw[color=gray] (-17,0) rectangle (-16,-6);
\filldraw[color=gray] (-15,0) rectangle (-14,-6);
\filldraw[color=gray] (-14,0) rectangle (-12,-5);
\filldraw[color=gray] (-11,0) rectangle (-9,-4);
\filldraw[color=gray](-9,0) rectangle (-8,-3);
\filldraw[color=gray] (-7,0) rectangle (-5,-2);
\filldraw[color=gray] (-5,0) rectangle (-4,-1);
\draw[thick]  (-17,0) grid (-4,-1);
\draw[thick]  (-17,-1) grid (-5,-2);
\draw[thick]  (-17,-2) grid (-7,-3);
\draw[thick]  (-17,-3) grid (-9,-4);
\draw[thick]  (-17,-4) grid (-12,-5);
\draw[thick]  (-17,-5) grid (-14,-6);

\filldraw[color=gray] (0,0) rectangle (1,-4);
\filldraw[color=gray] (2,0) rectangle (5,-4);
\filldraw[color=gray] (6,0) rectangle (7,-4);
\filldraw[color=gray] (7,0) rectangle (9,-3);
\filldraw[color=gray] (10,0) rectangle (13,-3);
\filldraw[color=gray] (14,0) rectangle (17,-1);
\draw[thick] (0,0) grid (17,-1);
\draw[thick] (0,-1) grid (14,-2);
\draw[thick] (0,-2) grid (13,-3);
\draw[thick] (0,-3) grid (7,-4);

\node at (-10,-5.5) {$\lambda'$};
\draw[ultra thick, ->] (4,-4.75) -- (4,-7.25);
\node at (6.25,-5.75) {$\psi(\phi(\lambda'))$};
\draw[ultra thick, ->] (-3.5,-2) -- (-1,-2);
\node at (-2.75,-3) {$\phi(\lambda')$};

\filldraw[color=gray] (1,-8) rectangle (14,-9);
\filldraw[color=gray] (1,-9) rectangle (11,-10);
\filldraw[color=gray] (1,-10) rectangle (11,-11);
\filldraw[color=gray] (1,-11) rectangle (6,-12);

\draw[thick] (1,-8) grid (14,-9);
\draw[thick] (1,-9) grid (11,-10);
\draw[thick] (1,-10) grid (11,-11);
\draw[thick] (1,-11) grid (6,-12);

\node at (.5,-8.5) {$1$};
\node at (.5,-9.5) {$2$};
\node at (.5,-10.5) {$1$};
\node at (.5,-11.5) {$3$};
\end{tikzpicture}
\caption{Shown here is the map $\lambda \mapsto \psi(\phi(\lambda'))$ applied to the partition (6, 6, 6, 5, 5, 4, 4, 4, 3, 3, 2, 2, 1), with $m=4$
and a  Schmidt weight that counts parts $\lambda_i$ with indices $i \not \equiv 2\; (\mathrm{mod}\; 4)$.}
\label{AndrewsKeith_main_example}
\end{figure}

With a fixed Schmidt weight, define the map $\psi(\lambda) = \mu$ on a partition $\lambda \in \mathcal{P}$ with Schmidt weight taken in conjugate, as follows.
For each $1 \leq i \leq \ell(\lambda)$, let $c_i$ be the remainder after dividing $\lambda_i$ by $m$. Form $\mu$ from $\lambda$
by removing all columns not counted by the Schmidt weight, and assign the part $\mu_i$ the color $c_i$. With this notation, the bijection given
in \cite{AndrewsKeith} between the two families of partitions in Theorem \ref{AndrewsKeith_main} can be written $\lambda \mapsto \psi(\phi(\lambda'))$ 
(see Figure \ref{AndrewsKeith_main_example} for an example of this map).

We now have the necessary tools to prove Theorem \ref{AndrewsKeith_extension} by extending this map. Before proceeding however, it is worth pointing
out that Theorem \ref{AndrewsKeith_main} actually admits a far more simple proof which was found by the author in the course of writing this paper, using
only elementary facts about generating functions and no difficult algebra. \newline

\noindent \textbf{A second proof of Theorem \ref{AndrewsKeith_main}.}
We start by proving the $m+1$ variable generating function identity
\begin{equation*}
\begin{aligned}
\smashoperator{\sum_{\lambda \in \mathcal{P}}} &z_1^{\displaystyle \rho_1} \cdots z_{m}^{\displaystyle\rho_{m}}
q^{\displaystyle \lambda_{s_1} + \cdots + \lambda_{s_i} + \lambda_{s_1+m} + \cdots + \lambda_{s_i+m} + \cdots} \\
&= \prod_{j=1}^m \left (q^{\lvert\{k \, \colon \, s_k \in S, \, k \leq j \}\rvert}z_j; q^i \right)^{-1}_\infty 
\end{aligned}
\end{equation*}
where $S =\{s_1, \dots, s_i\}$ and the $\rho_j$ are defined as in Theorem \ref{AndrewsKeith_main}, allowing in addition here $j=m$.
For a concrete example of this, let $m=5$ and $S=\{1, 3, 4\}$. Then this identitiy becomes
\begin{equation*}
\begin{aligned}
\smashoperator{\sum_{\lambda \in \mathcal{P}}} &z_1^{\displaystyle \rho_1} \cdots z_{5}^{\displaystyle\rho_{5}}
q^{\displaystyle \lambda_1 + \lambda_3 + \lambda_4 + \lambda_6 + \lambda_8 + \lambda_9 + \cdots} 
\\ &= \frac{1}{(qz_1;q^3)_\infty(qz_2; q^3)_\infty(q^2z_3; q^3)_\infty(q^3z_4; q^3)_\infty(q^3z_5; q^3)_\infty.}  
\end{aligned}
\end{equation*}

Any partition $\lambda \in \mathcal{P}$ is uniquely determined by independently choosing how many times each column height
appears in $\lambda$'s Young diagram.
Observe that each $\rho_j$ is counting the number of times columns in $\lambda$'s Young diagram with heights congruent to $j$ modulo $m$ appear.
Consider the partitions with only the column height $mn +j$ allowed, and let $K = \lvert\{k \, \colon \, s_k \in S, \, k \leq j \}\rvert$.
Clearly each column has Schmidt weight $in + K$, and so the generating function for these partitions is
\begin{equation*}
\frac{1}{1-z_jq^{ni +K}}.
\end{equation*}
Next, the generating function for all partitions with column heights all congruent to $j$ modulo $m$ can be obtained by taking the product over all $n \geq 0$,
which is
\begin{equation*}
\frac{1}{(z_j q^K; q^i)_\infty}.
\end{equation*}
Now take the product over all $j$ to get the identity.

Any partition $\lambda \in \mathcal{P}$ can be uniquely decomposed into a pair of partitions, with one in $\mathcal{D}_m$ and the other 
having parts repeating only in multiples of $m$. Simply combining the parts from the pair retrieves $\lambda$, and since the Schmidt weight
repeats with a period of $m$, any way this is done preserves the Schmidt weight. Hence we have the relation
\begin{equation}
\begin{aligned}
\smashoperator{\sum_{\lambda \in \mathcal{P}}} &z_1^{\displaystyle \rho_1} \cdots z_{m}^{\displaystyle\rho_{m}}
q^{\displaystyle \lambda_{s_1} + \cdots + \lambda_{s_i} + \lambda_{s_1+m} + \cdots + \lambda_{s_i+m} + \cdots} \\
&= \frac{1}{(z_mq^i; q^i)_{\infty}} \times \smashoperator{\sum_{\lambda \in \mathcal{D}_m}} z_1^{\displaystyle \rho_1} \cdots z_{m-1}^{\displaystyle\rho_{m-1}}
q^{\displaystyle \lambda_{s_1} + \cdots + \lambda_{s_i} + \lambda_{s_1+m} + \cdots + \lambda_{s_i+m} + \cdots}.
\end{aligned}
\label{Thm7_gf_2}
\end{equation}
This immediately implies Theorem \ref{AndrewsKeith_main} by canceling the common factor of $(z_mq^i;q^i)_\infty^{-1}$,
then interpreting the product in terms of colored partitions. \qed

\begin{remark}
The generating function identity established at the start of this proof has a partition identity interpretation that is similar to Theorem \ref{AndrewsKeith_main}, and
the map $\lambda \mapsto \psi(\lambda')$ gives a bijection between the two families of partitions. Theorem \ref{uncus_theorem_generalization} is also similar
ignoring the added statistics both identities carry. They become identical in the case $S=\{1\}$.
\end{remark}

\noindent \textbf{Proof of Theorem \ref{AndrewsKeith_extension}.}
Let $\lambda$ be as in the theorem statement. Define the pair of partitions $(\mu, \nu)$ as follows. For any part $\alpha_i$ appearing
$mk+r$ times for $0 \leq r < m$, add the $k$ parts $m\alpha_i$ to $\mu$, and $r$ parts $\alpha_i$ to $\nu$. Since we are moving identical parts in multiples of $m$,
the Schmidt weight of $\lambda$ is equal to that of $\mu'$ plus $\nu$. Observing also that $\nu \in \mathcal{D}_m$, we may apply $\phi$ to $\nu'$.
The bijection between the two families of partitions in the theorem is then $\lambda \mapsto \psi(\mu \cup \phi(\nu'))$ where $\mu \cup \phi(\nu')$
means combine the parts of $\mu$ and $\phi(\nu')$.

The fact that this map sends the Schmidt weight to the size of the colored partition is clear. Removing $m$ identical parts from $\lambda$ at a time
to form $\nu$ preserves all the $\rho_j$. Finally, $\psi$ only colors the parts of $\mu \cup \phi(\nu')$ in the color $m$ if they are multiples of $m$,
and only the parts from $\mu$ are by construction, so this map can be reversed. The remaining details are easily verified.
\qed

\section*{Acknowledgements}

The author would like to thank Professor William Keith for introducing him to this area of research, and for encouraging the author to present
this work in the Michigan Tech Specialty Seminar in Partition Theory, $q$-series, and Related Topics.


\begin{thebibliography}{1}\footnotesize

\bibitem{Alladi1} K. Alladi, Partition identities involving gaps and weights,
{\it Trans. Am. Math. Soc.} {\bf 349} (1997), 5001-5019.

\bibitem{Alladi2} K. Alladi, Schmidt-type theorems via weighted partition identities,
{\it Ramanujan J.} {\bf 61} (2023) 701-714.

\bibitem{Andrews} G. E. Andrews, \textit{The Theory of Partitions}, Cambridge
University Press, Cambridge, 1998.

\bibitem{AndrewsKeith} G. E. Andrews and W. J. Keith, A general class of
Schmidt theorems, {\it J. Number Theory} {\bf 247} (2023), 75-99.

\bibitem{AndrewsPaule} G. E. Andrews and P. Paule, MacMahon's partition
analysis XIII: Schmidt type partitions and modular forms,
{\it J. Number Theory} {\bf 234} (2022), 95-119.

\bibitem{BerkovichUncu} A. Berkovich and A. Uncu, On finite analogs of Schmidt's
problem and its variants, {\it S\'em. Lothar. Combin.} {\bf 88}, (2023).

\bibitem{Bessenrodt} C. Bessenrodt, A bijection for Lebesgue's partition
identity  in the spirit of Sylvester, {\it Discrete Math. } {\bf 123} (1994), 1-10.

\bibitem{Boulet} C. E. Boulet, A four-parameter partition identity, {\it Ramanujan J.} {\bf 12} (2006), 315-320.

\bibitem{Bowman} D. Bowman, Partitions with numbers in their gaps,
{\it Acta Arith.} {\bf 74} (1996), 97-105.

\bibitem{BridgesUncu} A. Bridges and A. Uncu, Weighted cylindric partitions,
{\it J. Algebraic Combin. } (2022).

\bibitem{ChernYee} S. Chern and A. J. Yee, Diagonal hooks and a Schmidt-type
partition identity, {\it Electron. J. Combin} {\bf 29} (2022).

\bibitem{DousseKim1} J. Dousse and B. Kim, An overpartition analog of $q$-binomial coefficients,
{\it Ramanujan J.} {\bf 42} (2017), 267-283.

\bibitem{DousseKim2} J. Dousse and B. Kim, An overpartition analogue of q-binomial coefficients, II: Combinatorial proofs and (q,t)-log concavity,
{\it J. Combin. Theory Ser. A} {\bf 158} (2018), 228-253.

\bibitem{Franklin} F. Franklin, On Partitions, {\it John Hopkins Univ. Cir.} {\bf 22} (1883), p. 72.

\bibitem{Glaisher} J. W. L. Glaisher, A theorem in partitions, {\it Messenger of Math.} {\bf 12} (1883), 158-170.

\bibitem{Ji} K. Q. Ji, A combinatorial proof of a Schmidt type theorem of Andrews,
and Paule, {\it Electron. J. Combin.} {\bf 29} (2022).

\bibitem{LiYee} R. Li and A. J. Yee, Schmidt type partitions,
{\it Enumer. comb. appl} {\bf 3} (2023).

\bibitem{CorteelLovejoy} S. Corteel and J. Lovejoy, Overpartitions
\textit{Trans. Amer. Math. Soc.} \textbf{356} (2004), Article no. 2.

\bibitem{Mork} P. Mork, Interrupted partitions, solution to problem 10629,
{\it Amer. Math. Monthly} {\bf 104} (1999), 87-88.

\bibitem{Schmidt} F. Schmidt, Interrupted partitions, {\it Am. Math. Mon.}
{\bf 104} (1999), 87-88.

\bibitem{Stockhofe} D. Stockhofe, Bijektive Abbildungen auf der menge der
partitionen einer naturlichen zahl, Ph. D thesis, Bayreuth. Math. Schr. 10
(1982), 1-59.

\bibitem{Uncu} A. K. Uncu, Weighted Rogers-Ramanujan partitions and Dyson's crank,
{\it Ramanujan J.} {\bf 46} (2018), 579-591.

\bibitem{WangXu} A. Y. Z. Wang and Z. Xu, The minimal excludant and Schmidt's
partition theorem, {\it Discrete Math.} {\bf 346} (2023).

\bibitem{WangZhangZhou} A. Y. Z. Wang, L. Zhang, and B. Zhou,
More on Schmidt's partition theorem, {\it Discrete Math.} {\bf 346} (2023).

\end{thebibliography}
\end{document}